\documentclass[11pt,twoside]{article}

\usepackage{color}
\usepackage{amsfonts}
\usepackage{amsmath}
\usepackage{amssymb}

\newcommand{\bean}{\begin{eqnarray*}}
\newcommand{\eean}{\end{eqnarray*}}
\newcommand{\benu}{\begin{enumerate}}
\newcommand{\eenu}{\end{enumerate}}
\newcommand{\eea}{\end{eqnarray}}
\newcommand{\bea}{\begin{eqnarray}}

\newtheorem{Theorem}{Theorem}

\newcommand{\be}{\begin{equation}}
\newcommand{\ee}{\end{equation}}

\newcommand{\N}{{\mathbb N}}

\newcommand{\R}{{\mathbb R}}

\newcommand{\SSS}{{\mathbb S}}

\newcommand{\la}{\lambda}

\newcommand{\e}{\epsilon}

\def\o{\circ}

\newcommand{\ben}{\begin{enumerate}}
\newcommand{\een}{\end{enumerate}}
\newcommand{\bit}{\begin{itemize}}
\newcommand{\eit}{\end{itemize}}
\newcommand{\edoc}{\end{document}}

\newcommand{\bdefi}{\begin{definition}}
\newcommand{\btheo}{\begin{theorem}}
\newcommand{\bprop}{\begin{proposition}}
\newcommand{\brema}{\begin{remark}}
\newcommand{\bcoro}{\begin{corollary}}
\newcommand{\blemm}{\begin{lemma}}
\newcommand{\bexam}{\begin{example}}

\newcommand{\edefi}{\end{definition}}
\newcommand{\etheo}{\end{theorem}}
\newcommand{\eprop}{\end{proposition}}
\newcommand{\erema}{\end{remark}}
\newcommand{\ecoro}{\end{corollary}}
\newcommand{\elemm}{\end{lemma}}
\newcommand{\eexam}{\end{example}}

\newcommand{\V}{\noindent}
\newcommand{\ci}{\circ}

\def\qed{\ensuremath{\quad\Box\quad}}

\title{On the non-slippery slope: Some observations on a recent paper on rolling bodies published in {\em Nature}}

\author{Olaf M\"uller\footnote{Humboldt-Universit\"at zu Berlin, Institut f\"ur Mathematik, Unter den Linden 6, 10099 Berlin. Email: mullerol@math.hu-berlin.de}}

\date{\today}

\begin{document}

\maketitle

\begin{abstract}
\V An interesting recent paper in {\em Nature} explores a new method to construct solid bodies rolling along given curves on an inclined plane, based on the Gauss Theorem. The present article complements those examinations with rigorous existence theorems and connections to seemingly unrelated questions like maritime rescue operations. 
\end{abstract}  

\section{Introduction and statement of the main results}

A recent article \cite{SDTEGG} published in {\em Nature} focuses on a very natural question: Given a curve $c$ on an inclined plane $P := A^{-1} (0), A(x) := x_1 + x_3 $, can we construct a solid body moving in time $t \mapsto K (t)$ (modelled as the support of a mass distribution for each $t$) that would, in the potential field $V : x \mapsto x_3 \forall x \in \R^3$ and the constraint $K(t) \subset A^{-1} ([0; \infty)) \ \forall t \in [0; \infty)$, move along $c$, i.e., such that its set of intersection with $P$ is approximately on $c$. As in \cite{SDTEGG}, we use a starshaped solid body: For a smooth "shape" function $S: \SSS^2 \rightarrow (0, \infty) $ define $P: \R^3 \setminus \{ 0|\} \rightarrow \SSS^2 , P(v) := v / \vert \vert v \vert \vert $ and   

\medskip

$K_S := \{ 0 \} \cup \{  v \in \R^3 \setminus \{ 0\} : \vert \vert v \vert \vert \in [0; S(P(v))]\}$

\medskip

\V The motion of a rigid body $B$ in an open subset $A$ of $\R^3$ can be described by a smooth curve $C:= (v,U): [0; S] \rightarrow D \subset \R^3 \times SO(3)$ where the first factor is to be understood as translations (acting on the center of mass), the second factor as the rotational degrees of freedom, and $D:= \{ (v,U)  \in \R^3 \times SO(3) | UB + v \subset A \}$, an open subset of $\R^3 \times SO(3)$.  

\bigskip

\V The equation obeyed by the rolling body is then

\bea
\label{MainEq}
C''(t) = C_0'' (t) - F C '(t)
\eea

\V where $C_0 '' (t) $ is the acceleration without rolling resistance obtained by $ E_{kin} '(t) = - E_{pot}'(t)$ (resulting from the Lagrange equation $\frac{d}{dt} \frac{\partial L}{\partial p} = \frac{\partial L}{\partial x} $ for the canonical local coordinates $x,p$ of the cotangent space) and $F = N \rho/r$ where $N$ is the normal component of the gravitational force, $\rho$ the rolling resistance coefficient, and $r$ the distance from the contact point to the center of mass. (What follows is equally valid for a broad class of rolling resistance models, as long as the rolling resistance is monotonously increasing in the velocity and collinear to $C'(t)$). $I_\alpha: \R^2 \rightarrow \R^3$ be defined by 

\medskip

$I_\alpha  :=   \left( \begin{array}{rrr}
\cos \alpha & 0 & \sin \alpha \\
0 & 1 & 0 \\
- \sin \alpha & 0 & \cos \alpha \\
\end{array} \right) 
\left( \begin{array}{rr}
	1&0\\
	0&1\\
	0&0\\
	\end{array} \right) =  \left( \begin{array}{rr}
	\cos (\alpha)&0\\
	0&1\\
	- \sin (\alpha)&0\\
	\end{array} \right),$

\V and let $J: T(I_\alpha (\R^2)) \rightarrow T(I_\alpha (\R^2))$ be defined by $J( x) := I_\alpha (i (I_\alpha^{-1}) (x))$ for all $x \in I_\alpha (\R^2)$, where we use an identification of $\R^2 $ with $\mathbb{C}$.
.
\medskip  

\begin{Theorem}
	\label{Main}
Let $f \in C^2 (\R , \R)$ be $T$-periodic and let $\epsilon>0$. Then there is a ball radius $R \in (0; \infty)$, a rolling resistance coefficient $C>0$ and a slope angle $ \alpha \in (0; \pi/2)$, which can all three be estimated from above in terms of $\vert \vert f \vert \vert_{C^2}$, and a smooth function $S: \SSS^2 \rightarrow (1 - \e ; 1]$, such that the solution $F_{\alpha}: [0, \infty) \rightarrow  {\rm Isom} (\R^3)$ to Equation \ref{MainEq} for the characteristic function $\chi_{K_{RS}}$ of the set $K_{RS}$ satisfies after appropriate reparametrization by some increasing $ u \in C^1 (\R)$ that there is $ \delta \in (0; \epsilon)$ with

$$F_{\alpha} (\chi_{K_{RS}}) (u(t))  \cap I_\alpha (\R^2) = I_\alpha ((t,f(t))) + [- \delta; \delta ] \cdot J(c'(t)) \ \ \forall t \in \R.$$

\end{Theorem}

\V{\bf Remark.} The conclusion of the theorem means that the contact set between the rolling body $K_{RS}$ and the inclined plane $I_\alpha (\R^2)$ at time $t$ is a segment centered at $ I_\alpha (c(t)) $ of width $2 \delta$, for all $t$. The map $\tau$ assigning to each $S $ the curve of the barycenters\footnote{The bounded set $ F_{\alpha} (u(t)) (K_S) \cap I_\alpha (\R^2)$ is also a closed subset of $\R^3$ and therefore Lebesgue measurable.} of the contact sets of $K_S$ at the times $t$ is continuous w.r.t. the $C^2$ topology on the left-hand side and the compact-open topology ($C^0$ convergence on each compact subinterval) on the left-hand side, and the statement of the theorem provides a continuous right inverse of $\tau$. However, $\tau$ is not injective, as shown by the simple example of sphere and cylinder both having an affine trajectory. Mind that the rolling body's mass distribution is taken here to be homogeneous on $K_{RS}$ as opposed to the setting of \cite{SDTEGG} where the mere existence of such a mass distribution was shown, typically with the mass concentrated in a heavy core.

\V {\bf Remark.} Without taking into account rolling resistance, the statement would be wrong for every curve except a line: By energy preservation, the  body's velocity would tend to infinity with time, and then centrifugal forces would make the body's leave its course. {\em With} rolling resistance, however, we have to assume that the curve is always downhill: Not only that some courses are impossible to realize due to energy preservation, we also have strange discontinuities: Consider the case of a cone, then its course depends on whether it turns around the tip or around the opposite boundary, and this dichotomy depends  non-continuously on the rolling resistance. This existence of tipping points makes the entire analysis quite complicated. 

\bigskip	

\V We reduce, like in \cite{SDTEGG}, the problem to the finding a simple\footnote{Simpleness of the curve is not considered in \cite{SDTEGG}, however needed for the "carving".} closed $r$-lifting $c_r := L_r(c)$ of $c$. The $r$-lifting $L_r(c): \R\rightarrow \SSS^2 (r)$ is defined by $L_r(c) (0) = (r,0,0)$, $(L_r(c))'(0) \in \R \cdot (0,1,0)$, $|| (L_r (c) )'(t) || = ||c'(t||)$ for all $t$ and $\kappa_c (t) = \kappa_{L_r (c) (t)} $ for all $t$, where $\kappa$ is the curvature. This encodes rolling without slipping or pivoting. Once we found such an $r$-lift, fixed in the following, we find a normal neighborhood $N$ of width $b$ around it, then we define $S(x):= \min \{ b /r - COS(d(x, c(\R)) ) \}$ where $COS$ is the real function coinciding with $\cos$ on the positive numbers and being equal to $1$ on the nonpositive numbers. This yields a continuous function corresponding to the carving of a small flat trajectory around $L_r(c)$. As observed in the original article, a tilting movement (leaving the carved path) would lead to a cycloidal movement of the center of mass upwards, not differentiable at the origin, and thus not realized, as connected with an increase of energy, {\em as long as the center of mass is contained in a wedge of opening angle $\beta := \alpha \cdot a$ over the contact segment}. Here $a$ is an a priori $C^1$ bound of the function $f$. The volume of a spherical cap for a height $h$ is $v(h) = \frac{\pi h^2}{3} (3r-h)$, thus a very generous but sufficient estimate of the volume $V$ of the carved part is $V < v \cdot l$ with $l := \ell (c)$. An equally generous estimate of the distance of the barycenter of the resulting solid body for carving heights $0$ and $h$ is $|| B(h) - B(0) || < V \cdot r $, so all in all we get

$$ || \frac{B(h) - B(0)}{ h } || < \frac{\pi h}{3} (3r-h) \cdot r l < \pi h r^2 l  \rightarrow_{h \rightarrow 0} 0 .$$

\V As the wedge width $ b $ satisfies $b= \sqrt{h(h-r)}$, there is a small height $h$ such that the barycenter is contained in the respective wedge. This will allow to choose the mass distribution of the solid body homogeneous once we found a good corresponding curve on the sphere.
  
\bigskip

\V The central observation of \cite{SDTEGG} at this point is the following: Let $k: [0;T] \rightarrow \SSS^2$ be $C^0$ and piecewise $C^1$. If there is a point $p \in \SSS^2$ with $d(p,c(0)) = d(p, c(T))$ and the enclosed area $A(g_{c(T) p} * c * g_{p} * g_{p c(0)}) = \pi/n  $, then $\bigstar_{j=0}^{n-1} R^j \o k $ closed. The proof in \cite{SDTEGG} makes use of Gauß's Theorem (for a constant integrand). One could parallel this in a first step when we do not care about simple-closedness. Then we would have to (straightforwardly) generalize Gauß's Theorem to smooth maps $\phi$ from the unit disc to a surface using degree theory: The integrand in the image is then the (i.g. non-continuous) degree of $\phi$. However, in the following we do not use Gau\ss's Theorem. 

\bigskip

\V Now we recall that the carving process {\em requires that the $L_r(c)$ be not only closed but also simple (i.e., injective)}. Including simpleness in our requirements makes the entire procedure unstable in the sense that if we have found $r >0$ with $L_r (c)$ simple closed, we can modify $c$ by an arbitrarily small amount in any $C^k$ metric to obtain $\tilde{c}$ with such that appropriate restrictions of $L_s(\tilde{c})$ intersect each other transversally for each $s$ in an open interval around $r$.

\begin{Theorem}
	Let $k: [0;1] \rightarrow \R^2$ be $C^1$ and piecewise $C^2$, with $k'(0) = k'(1)$. Then there is $r>0$ such that $L_r (k)$ is simple and closed. 
\end{Theorem}

\V{\bf Proof.} Consider the concatenation $c$ of infinitely many translates of $k$ by $v:= k(1) - k(0)$. It is a $C^1$ and piecewise $C^2$ curve again with $c \o T_v (\R) = c (\R)$. Let $ b_r:= \min \{ d(c_r (s)), c_r (s+T) | s \in [0;T] \} $. Then $b_r$ depends continuously on $r$. By the spherical ASA congruence theorem, the angle $\alpha_r$ between the left normals $n,n'$ is determined (continuously) by $(r ,b_r) \in [0; \infty] \times (0; T]$. In the same time, the lengths of $n$ and $n'$ agree (both agree with the radius of the largest ball in the complement of smaller connected component in the complement of $c_r(\R)$ invariant under the symmetry of $c$). We want $\alpha_r \in  \pi/\N^* $. As have $\alpha_\infty = 0$ and $\alpha_r >0$ for all large $r < \infty$, the statement follows from the intermediate value theorem (for the simplicity statement use that $L_r - (r,0,0)\rightarrow_{r \rightarrow \infty} {\rm id}$). \hfill \qed

\bigskip

\V However, we cannot restrict $r$ in terms of $\ell(c)$:
 
\begin{Theorem}
Let $D>0$. Then there is $\gamma: [0;2 \pi] \rightarrow \R^2$ smooth, parametrized by arclength, such that $L_r (\gamma)$ is not injective for each $r <D$ and appropriate restrictions of $L_r(\gamma)$ intersect each other transversally.
\end{Theorem}

 \V {\bf Proof.} Consider an arclength-parametrized circle curve $c: t \mapsto (\cos (t), \sin (t))$ of radius $R:= 1$ in $\R^2$. Its curvature is $1/R = 1 $, its length $2 \pi $. Now consider curves $c_{r, \theta}$ of constant latitude $\theta$ in $\partial B(0, r) \subset \R^3$: $c_{r, \theta} (t) := r (\cos \theta \sin t , \cos \theta \cos t, \sin \theta) \ \forall t \in [0; 2 \pi]$. Their curvature is easily seen to be constantly $k_{r, \theta} := \frac{1}{r} \tan \theta$, their length is $\ell_{r, \theta} := 2 \pi r \cos \theta$. This implies that locally $L_r (c)$ agrees with $c_{r, \theta(r,R)}$ where $\theta (r, R) := \arctan (r/R)$. Thus,

$$\ell (L_r (c_R))  = 2 \pi r \cos (\arctan (r/R)) = 2 \pi \frac{1}{\sqrt{\frac{1}{r^2} + \frac{1}{R^2} }} < 2 \pi R \ \  \forall r >0 ,$$ 

\V and we easily see that there is even an infinite-dimensional subspace $V$ of $C^2 ([0;2\pi], \R^2)$ and $\e >0$ such that for all $\gamma$ in $c_R + (B (0,\e) \cap V)$ we obtain that $L_r (\gamma) $ is not injective, and appropriate restrictions of $L_r(\gamma)$ intersect each other transversally. \hfill \qed 

\bigskip

\V As argued above, it is reasonable to restrict to those curves $c$ whose images are functions of the $x_1$-coordinate. Define, for $f \in C^0 ( [0;\la], \R) $, $c_f := A(f): [0; T(f)] \rightarrow \R^2$ as the injective arclength-parametrized curve with $ c_f (0) = (0, f(0)) $ and $c_f ([0; T(f)]) = f \in \R^2$, in particular $T(f) = \ell (c_f)$. For each $r>0$ we define the curves  $c_{f,r} := L_r (c_f): [0; T] \rightarrow \SSS^2(r)$ with $c_{f,r}(0) = (r,0,0) $, $c_{f,r}'(0) = (0,r,0)$, $|| c_{f,r}'(t) || = 1 \forall t \in [0; T(f)]$ and $\kappa_{c_{f,r}} (t) = \kappa_{c_f}(t) \forall t \in [0; T(f)]$ (equality of the respective curvatures); we also define $F_{r}: C^0 ([0;T]) \rightarrow \R$ by $ F_r(k) = d(C_{\kappa,r}(0) , c_{k,r} (T))$. First, even in this case we need to impose the following linear bound of $r$ in terms of $\ell$ if we insist on injectivity:

\begin{Theorem}
\label{Gegenbsp-inj}	
Let $a := \pi \cdot \sqrt{\frac{\sqrt{17} -1}{2}} $. For each $T>0$ there is a smooth function $f$ of $\ell(c_f)= T$ such that for each $r < T/a$,  the curve $c_{f,r}$ is not injective.  	
\end{Theorem}	

\V {\bf Remark.} As $a \in (3.9;4)$, we get e.g. the more memory-friendly estimate that in general we need $r > \ell (c_f)/4$ to ensure injectivity of $c_{f,r}$.

\V{\bf Proof.} For $f: x \mapsto \sqrt{\rho^2 - x^2}$ on $[- \rho; \rho]$, $c_f$ is a semi-circle of radius $\rho$, arc length $\ell = \pi \rho$ and constant curvature $\kappa = \rho^{-1} = \pi/ \ell$. The corresponding arclength-parametrized curves $c_{(\ell)} $ (which are the $c^f$ of above) of total length $\ell$ are $ c_{(\ell)} : t \mapsto \ell/\pi \cdot (\sin (t \pi/\ell), \cos(t \pi/\ell ))$. The curve $c_{\ell, r} := (c_{(\ell)})_r $ are the circles of latitude $\theta$ of the same constant curvature in $\partial B(0,r)$. The circle of latitude $\theta $ in $\partial B(0,r)$ has curvature $ \kappa_r := 1/(r \tan \theta)$ and length $\ell_r := 2 \pi r \cos \theta = 2 \pi r \cos \arctan (\frac{\ell}{\pi  r}) = \frac{2 \pi r }{\sqrt{1 + \frac{\ell^2}{\pi^2 r^2}}} = :\sigma (r)$. The function $r \mapsto \sigma (r)$ is increasing, and $\sigma (r) = \ell$ is, with $a:= \frac{\ell^2}{\pi^2 r^2}$, equivalent to the quadratic equation $a^2 +a - 4 = 0$, which has the roots $a_\pm := \frac{-1 \pm \sqrt{17}}{2}$, and by definition of $a$, the root $a_+$ is the one relevant to our purpose. The rest is reinsertion of definitions. \hfill \qed

\begin{Theorem}[Man-over-board curve]
	\label{ManOverBoard}
Let $D: C^1 ([0, \la]) \rightarrow C^0([0; \la])$ be the operator of differentiation $D(f) = f'$. If $\la < \pi r $, then $F_r$ achieves its minimum on $\kappa^{-1}(D(B_{L^\infty}(0,R) )) \cap C^2([0; \la], \R^2)$ at the constant function $\la/R$.  	
\end{Theorem}

\V {\bf Remark.} The bound on $\la$ entails that the image of $c_r$ is contained in the hemisphere around $c(\la/2)$. The result generalizes the well-known (but not well-documented) statements that $F_\infty$ achieves its minimum on $B_{L^\infty} (0,R/\la) \cap C^2 ([0; \la], \R^2)  $ (acceleration bound) and also on $B_{L^1} (0,R/\la) \cap C^2 ([0; \la], \R^2)$ (fuel bound) at the constant function $\la /R$, which is put in practice in man-over-board manoeuvers. But $D(B_{L^\infty}(0, R)) \subsetneq B_{L^1} (0, R \cdot \la)$, thus the statement of the theorem above is more general then the two folk theorems. 

\V {\bf Proof.} We first write the parageodesic equation (the equation for a arclength-parametrized curve of prescribed curvature) in spherical coordinates $(\theta, \phi)$ around some point $c(0)$, i.e., with $\theta (c(0)) = 0$, which together with the definition of $\theta, \phi$ yields $\det (c'(t), \partial \phi (c(t))) >0$ for small nonzero $t$. The problem can easily be reduced to the one of minimization of the first time $T$ after which $\langle c' , \partial \phi \rangle = 0  $ in the respective class of curves (If we solved the problem of $T$-mimimizing by a curve $k$, then the original minimizing problem is solved by gluing $k$ to $R \o k$ where $R$ is the reflection at the line containing the endpoints of $k$). We use the Christoffel symbols 

$$ \Gamma^\theta = \left( \begin{array}{rr}
	0&0\\
	0& - \sin \theta \cos \theta \\
\end{array} \right)  ,   \Gamma^\phi = \left( \begin{array}{rr}
0& \cot \theta \\
\cot \theta & 0 \\
\end{array} \right) ,  $$

\V for the parageodesic equation for $c = (\theta, \phi)$, which reads 

$$ \left( \begin{array}{r}
	\theta'' (t) \\
	 \phi '' (t))\\
\end{array}   \right)  + \left( \begin{array}{r}
\Gamma^\theta (c'(t), c'(t)) \\
\Gamma^\theta (c'(t), c'(t))\\
\end{array}   \right)  = \nabla_t c'(t) = \kappa (t) \cdot \left( \begin{array}{rr}
0 & \sin (\theta) \\
1/\sin(\theta) & 0 \\
\end{array}   \right) \cdot c'(t)  $$

\V to see that in any interval $J:= [0;T]$ in which $\theta' >0$, the curve $\tilde{c} := c_{|k|}  $ satisfies $\theta \o \tilde{c} |_J \leq \theta \o c |_c $, (thereby reducing the $T$-minimization problem to the fuel-bound problem). This implies via the first integral of arc-length:

$$ \theta '' (t) = - \kappa (t) \sin (\theta) \sqrt{1 - (\theta ')^2} - \sin (\theta) \cos (\theta) (1 -  (\theta')^2) , $$

\V and the classical result about the differential dependence of an ODE solution on a control parameter (here $\kappa$) (see e.g. \cite{amann}, Th. II.9.2), tell us that the derivative of a solution along a positive variation of $\kappa$ is positive, so a constant function $\kappa$ is actually the minimum in the respective spaces. \hfill \qed

\bigskip 

 \V Then we get a statement of injectivity of $L_r(c)$ if we impose additionally that the curve is function-like, i.e., the image of $c$ is a function:

 \begin{Theorem}
 	\label{MOB-Corollary}
 	Let $f \in C^2([0;E])$, then for each $r > \ell (c_f)$, the curve $c_{f,r} $ is injective. Moreover, for $r \geq \sqrt{2} \ell(c_f)$, we get $F_r > \ell(c_f)/3 $.
 \end{Theorem}

\V{\bf Proof.} First, for the arc-length parametrization $A$ as defined before Theorem \ref{Gegenbsp-inj}, every $c \in A(C^2 ([0 \la]))$, i.e., each curve $c_f$ for some function $f$, satisfies $\kappa_c \in D(B(0, \pi))$ by applying the fundamental theorem of differential calculus to $\langle c', e_2 \rangle$. The result follows by comparison of the semicircle $c_f$ of radius $\rho$ with curvature $k_\infty = \rho^{-1}$ and length $\ell_{\infty} = \pi \rho$ with the curves of latitude $\theta $ in $\partial B(0,r)$ with curvature $K_r = \frac{1}{r \tan \theta}$ and length $\Lambda_r = 2 \pi r \sin \theta$. Equating the curvatures yields $2 \pi r \sin \theta  = 2 \pi \rho / \sqrt{1 + \frac{\rho^2}{r^2}}$. Now following Theorem \ref{ManOverBoard}, injectivity of a curve $c$ with $\kappa_c \in D(B(0, \pi))$ is satisfied if $\Lambda_r > \ell_\infty$, which is equivalent to $r > \rho$. This is the case if $r > \ell_\infty (c) $ (in which case $r > \ell_{\infty} (c) / \pi$, which was one of the conditions). For the second assertion, we replace the requirement $\Lambda_r > \ell_\infty$ with $\frac{3}{4}\Lambda_r > \ell_\infty$ to see that $r \geq \sqrt{2} \rho$ implies that $F_r$ is smaller or equal to $ G_r  $, the distance of the endpoint of three quarters of a curve of constant latitude and of length $\ell_\infty$ in the sphere of radius $r$. This distance can be strictly estimated above by the length of the missing quarters, which yields a length ratio of $1/3$. \hfill \qed

\bigskip

\V The last result is that with a $C^2$ bound on $c_\infty$ we can indeed find a bound on the diameter of a rolling body following $c_\infty$. The same is true with a bound of the form $\kappa \in D(B(0,R))$ as in Th. \ref{ManOverBoard}, allowing the curve to be arbitrarily wildly curved but only locally, without collecting too much curvature. However, in practice a bound on the curvature is needed anyway in order to get a bound on the friction needed to keep the rolling body in the trace.

\begin{Theorem}
	\label{SimpleClosed}
	For each $D>0$ there is $\rho >0$ such that for all $ f \in C^2 ([0;T])$ with $ f'(0) = 0 = f'(T)$ and $||f||_{C^2} <D$ there is an $r \in (0; \rho)$ such that $L_r (c_f)$ is simple closed.  
\end{Theorem}

\V{\bf Proof.} As before, we loosely identify $c$ and its periodic continuation. The $C^2$ bound entails $ |\kappa (t)| \leq D $ for all $t \in \R$ and yields a bound from below on 

\bean
b_r : = \min \{ d(p, c_r(t)) : t \in \R, d(p, c_r(t) ) = d(p, c_r(t + T)) \\ 
\land  \ p - c_r(t)  \perp c_r'(t) \ \land \ p - c_r(t + T)  \perp c_r'(t+T)  \}
\eean 

\V In fact, as $c_r$ osculates in $c(t)$ the circle of radius $b_r$ around $p$, which has curvature $\frac{1}{r \tan (b_r/r)}$, we get $D^{-1} < r \tan (b_r/r) $, i.e., $b_r > r \arctan (1/rD)$. Moreover, $F_r := d(c_r(t), c_r(t+T)) \leq \ell (c_\infty) \leq T \sqrt{1+D^2}$ does not depend on $t$. Also, the angle $a_r:= \angle (c(t) p c(t+T))$ does not depend on the choice of such a point with minimal $b_r$, by the spherical SSS congruence theorem. By the spherical law of cosines,

\[ \cos (a_r) =  \frac{\cos (F_r/r) - \cos^2 ( b_r/r)}{\sin^2 (  b_r/r)}. \]

 \V Theorem \ref{MOB-Corollary} also yields a lower bound for $F_r$, an upper bound for both $F_r$ and $b_r$ is given by $\pi r/2$, and a lower bound to $ b_r $ has been calculated above. This allows to increase the angle at $p$ by at least a factor of $3$ if we decrease $r$ by a factor of $6$. Then we use 

\[\forall a \in [0; \pi] \exists n \in \N^*:  2 \pi /n  \in [a/3;a]   \]

\V to show that $\bigstar_{k=1}^n R^k \ci c$ is closed for a rotation $R$ around $p$. \hfill \qed 

\bigskip

\V Theorem \ref{Main} now follows from Theorem \ref{SimpleClosed} and the arguments before. \hfill \qed

\end{document}